\documentclass[12pt]{amsart}
\usepackage{epsfig}

\DeclareFontEncoding{OT2}{}{} 
\newcommand{\textcyr}[1]{%
 {\fontencoding{OT2}\fontfamily{cmr}\fontseries{m}\fontshape{n}\selectfont #1}}

\newcommand{\Sha}{{\mbox{\textcyr{Sh}}}}

\newcommand{\Q}{{\mathbf Q}}
\newcommand{\Qbar}{{\overline{\Q}}}

\newcommand{\Fbar}{{\overline{\F}}}

\newcommand{\Wbar}{{\overline{W}}}
\newcommand{\Z}{{\mathbf Z}}
\newcommand{\R}{{\mathbf R}}
\newcommand{\Aff}{{\mathbf A}}
\newcommand{\PP}{{\mathbf P}}

\newcommand{\F}{{\mathbf F}}

\newcommand{\sw}{\operatorname{sw}}
\newcommand{\sing}{{\operatorname{sing}}}

\newcommand{\Spec}{\operatorname{Spec}}

\newcommand{\HH}{{\operatorname{H}}}
\newcommand{\RR}{{\operatorname{R}}}
\newcommand{\et}{{\operatorname{\acute{e}t}}}
\newcommand{\isom}{\cong}
\newcommand{\del}{\partial}

\newcommand{\notdiv}{{\not |}}

\newtheorem{theorem}{Theorem}
\newtheorem{lemma}[theorem]{Lemma}

\newtheorem{prop}[theorem]{Proposition}

\theoremstyle{definition}

\theoremstyle{remark}

\begin{document}

\title[Family of genus-one curves]{An explicit algebraic family of genus-one curves violating the Hasse principle}
\author{Bjorn Poonen}
\thanks{This research was supported by
National Science Foundation grant DMS-9801104,
an Alfred P. Sloan Fellowship, and a David and Lucile Packard Fellowship.}
\address{Department of Mathematics, University of California, Berkeley, CA 94720-3840, USA}
\email{poonen@math.berkeley.edu}
\date{October 19, 1999}

\begin{abstract}
We prove that for any $t \in \Q$,
the curve
	$$5 x^3 + 9 y^3 + 10 z^3
	+ 12 \left( \frac{t^{12}-t^4-1}{t^{12}-t^8-1} \right)^3
			(x+y+z)^3 = 0$$
in $\PP^2$ is a genus~1 curve violating the Hasse principle.
An explicit Weierstrass model for its Jacobian $E_t$ is given.
The Shafarevich-Tate group of each $E_t$ contains
a subgroup isomorphic to $\Z/3 \times \Z/3$.
\end{abstract}

\maketitle

\section{Introduction}
\label{introduction}

One says that a variety $X$ over $\Q$ {\em violates the Hasse principle}
if $X(\Q_v) \not= \emptyset$ for all completions $\Q_v$ of $\Q$
(i.e., $\R$ and $\Q_p$ for all primes $p$)
but $X(\Q)=\emptyset$.
Hasse proved that degree~2 hypersurfaces in $\PP^n$
{\em satisfy} the Hasse principle.
In particular, if $X$ is a genus~0 curve,
then $X$ satisfies the Hasse principle,
since the anticanonical embedding of $X$ is a conic in $\PP^2$.

Around 1940,
Lind~\cite{lind} and (independently, but shortly later)
Reichardt~\cite{reichardt} discovered examples of genus~1 curves over $\Q$
that violate the Hasse principle,
such as the nonsingular projective model of the affine curve
	$$2y^2 = 1-17x^4.$$
Later, Selmer~\cite{selmer} gave examples of diagonal plane cubic curves
(also of genus~1) violating the Hasse principle, including
	$$3 x^3 + 4 y^3 + 5 z^3 = 0$$
in $\PP^2$.

O'Neil in~\cite[\S6.5]{oneil} constructs an interesting example
of an algebraic {\em family} of genus~1 curves
each having $\Q_p$-points for all $p \le \infty$.
Some fibers in her family violate the Hasse principle,
by failing to have a $\Q$-point.
In other words,
these fibers represent nonzero elements of the Shafarevich-Tate groups
of their Jacobians.

In~\cite{colliotpoonen}, Colliot-Th\'el\`ene and the present author
prove, among other things,
the existence of non-isotrivial families of genus~1 curves
over the base $\PP^1$, smooth over a dense open subset,
such that the fiber over {\em each} rational point of $\PP^1$
is a smooth plane cubic violating the Hasse principle.
In more concrete terms, this implies that
there exists a family of plane cubics depending
on a parameter $t$,
such that the $j$-invariant is a non-constant function of $t$,
and such that substituting {\em any} rational number for $t$
results in a smooth plane cubic over $\Q$ violating the Hasse principle.

The purpose of this paper is to produce an explicit example of
such a family.
Our example, presented as a family of cubic curves in $\PP^2$ with
homogeneous coordinates $x,y,z$, is
	$$5 x^3 + 9 y^3 + 10 z^3
	+ 12 \left( \frac{t^{12}-t^4-1}{t^{12}-t^8-1} \right)^3
			(x+y+z)^3 = 0.$$

\section{The cubic surface construction}
\label{cubicsurface}

Let us review briefly the construction in~\cite{colliotpoonen}.
Swinnerton-Dyer~\cite{swinnertondyer}
proved that there exists a smooth cubic surface $V$ in $\PP^3$
over $\Q$ violating the Hasse principle; choose one.
If $L$ is a line in $\PP^3$ meeting
$V$ in exactly 3 geometric points,
and $W$ denotes the blowup of $V$ along $V \cap L$,
then projection from $L$ induces a fibration $W \rightarrow \PP^1$
whose fibers are hyperplane sections of $V$.
Moreover, if $L$ is sufficiently general,
then $W \rightarrow \PP^1$ will be a Lefschetz pencil,
meaning that the only singularities of fibers are nodes.
In fact, for most $L$,
all fibers will be either smooth plane cubic curves,
or cubic curves with a single node.

For some $N \ge 1$,
the above construction can be done with models over $\Spec \Z[1/N]$
so that for each prime $p \notdiv N$, reduction mod $p$
yields a family of plane cubic curves
each smooth or with a single node.
One then proves that if $p \notdiv N$,
each fiber above an $\F_p$-point has a smooth $\F_p$-point,
so Hensel's Lemma constructs a $\Q_p$-point
on the fiber $W_t$ of $W \rightarrow \PP^1$ above any $t \in \PP^1(\Q)$.

There is no reason that such $W_t$ should have $\Q_p$-points
for $p|N$, but the existence of $\Q_p$-points on $V$
implies that at least for $t$ in a nonempty $p$-adically open subset $U_p$
of $\PP^1(\Q_p)$, $W_t(\Q_p)$ will be nonempty.
We obtain the desired family by
base-extending $W \rightarrow \PP^1$
by a rational function $f: \PP^1 \rightarrow \PP^1$
such that $f(\PP^1(\Q_p)) \subseteq U_p$ for each $p|N$.

More details of this construction can be found in~\cite{colliotpoonen}.

\section{Lemmas}
\label{lemmas}

\begin{figure}
\begin{center}
\epsfig{file=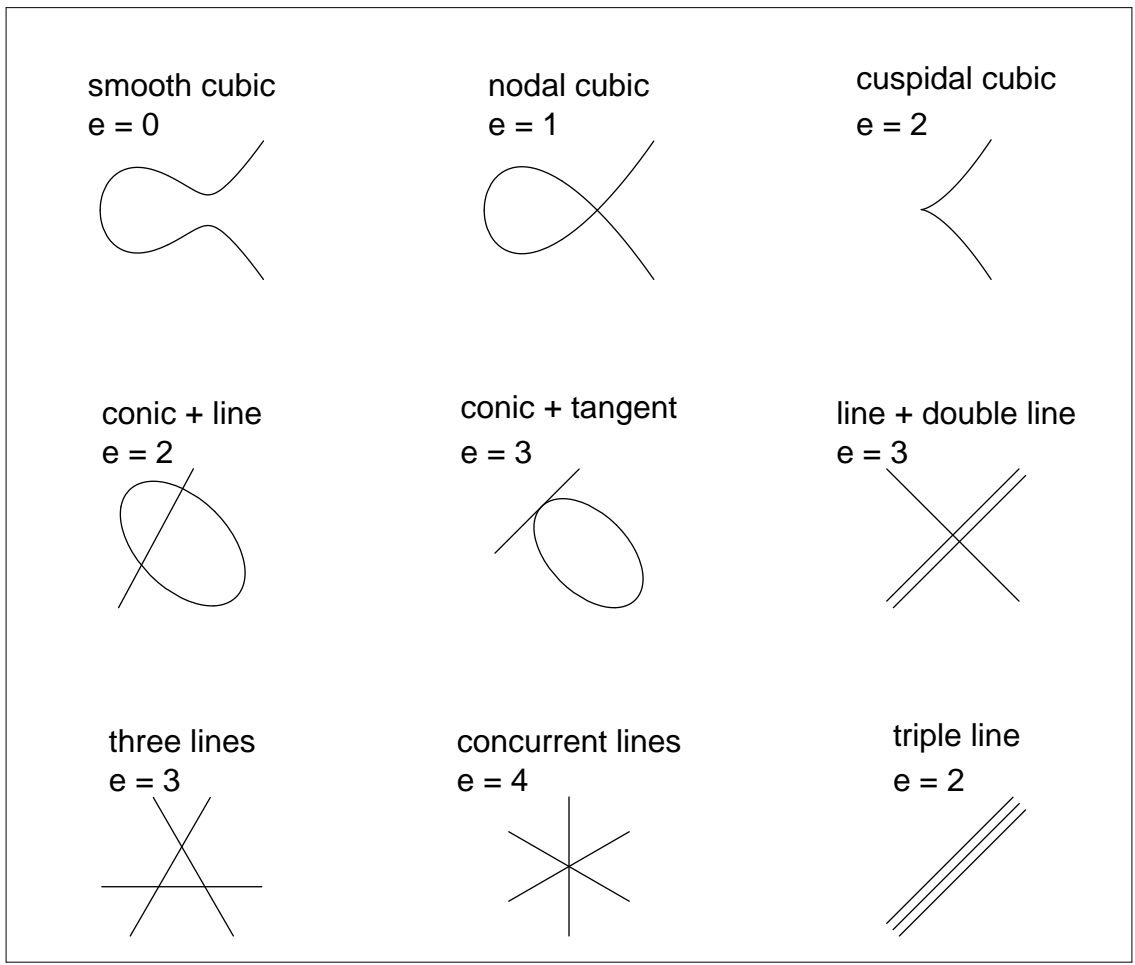} 
\end{center}
\caption{Plane cubic curves with their Euler characteristics.}
\label{planecubics}
\end{figure}

\begin{lemma}
\label{twelve}
Let $V$ be a smooth cubic surface in $\PP^3$
over an algebraically closed field $k$.
Let $L$ be a line in $\PP^3$ intersecting $V$ in exactly 3 points.
Let $W$ be the blowup of $V$ at these points.
Let $W \rightarrow \PP^1$ be the fibration
of $W$ by plane cubics
induced by the projection $\PP^3 \setminus L \rightarrow \PP^1$ from $L$.
Assume that some fiber of $\pi: W \rightarrow \PP^1$ is smooth.
Then at most 12 fibers are singular,
and if there are exactly 12,
each is a nodal plane cubic.
\end{lemma}

\begin{proof}
Let $p$ be the characteristic of $k$,
and choose a prime $\ell \not= p$.
Let
	$$\chi(V) = \sum_{i=0}^{2 \dim V}
			(-1)^i \dim_{\F_\ell} \HH_\et^i(V,\F_\ell)$$
denote the Euler characteristic.
Since $V$ is isomorphic to the blowup of $\PP^2$ at 6 points
(\cite[Theorem~24.4]{manincubic}, for example),
	$$\chi(V) = \chi(\PP^2) + 6 = 3 + 6 = 9.$$
Since $W$ is the blowup of $V$ at 3 points,
	$$\chi(W) = \chi(V) + 3 = 9 + 3 = 12.$$

On the other hand,
combining the Leray spectral sequence
	$$\HH^p(\PP^1,\RR^q \pi_\ast \F_\ell)
		\implies \HH^{p+q}(W,\F_\ell)$$
with the Grothendieck-Ogg-Shafarevich formula
(\cite[Th\'eor\`eme~1]{raynaud} or~\cite[Theorem~2.12]{milneetale})
yields
\begin{equation}
\label{ecformula}
	\chi(W) = \chi(W_\eta) \; \chi(\PP^1) + \sum_{t \in \PP^1(k)}
		\left[ \chi(W_t) - \chi(W_\eta)
			- \sw_t(\HH_\et^\ast(W_\eta,\F_\ell)) \right],
\end{equation}
where $W_\eta$ is the generic fiber,
$W_t$ is the fiber above $t$,
and
	$$\sw_t(\HH_\et^\ast(W_\eta,\F_\ell)):=
		\sum_{i=0}^2 (-1)^i \sw_t(\HH_\et^i(W_\eta,\F_\ell))$$
is the alternating sum of the Swan conductors of $\HH_\et^i(W_\eta,\F_\ell)$
considered as a representation of
the inertia group at $t$ of the base $\PP^1$.
Since $W_\eta$ is a smooth curve of genus $g=1$, $\chi(W_\eta)=2-2g=0$.
If $t \in \PP^1(k)$ is such that $W_t$ is smooth,
then all terms within the brackets
on the right side of~(\ref{ecformula}) are $0$,
so the sum is finite.
The Swan conductor of
	$$\HH_\et^0(W_\eta,\F_\ell) \isom \HH_\et^2(W_\eta,\F_\ell)
		\isom \F_\ell$$
is trivial.
Hence~(\ref{ecformula}) becomes
	$$12 = \sum_{t:\; W_t \text{ is singular}}
			\left[ \chi(W_t)
				+ \sw_t(\HH_\et^1(W_\eta,\F_\ell)) \right],$$
Since $\sw_t(\HH_\et^1(W_\eta,\F_\ell))$ is a dimension, it is nonnegative,
so the lemma will follow from the following claim:
if $W_t$ is singular, $\chi(W_t) \ge 1$ with equality
if and only if $W_t$ is a nodal cubic.
To prove this, we enumerate the combinatorial
possibilities for a plane cubic,
corresponding to the degrees of the factors of the cubic polynomial:
see Figure~\ref{planecubics}.
The Euler characteristic for each, which is unchanged if we pass to the
associated reduced scheme $C$, is computed using the formula
\begin{equation}
\label{normalization}
	\chi(C) = \sum_i (2-2g_{\tilde{C}_i})
		+ \# C_\sing - \# \alpha^{-1}(C_\sing),
\end{equation}
where $\alpha: \tilde{C} \rightarrow C$ is the normalization of $C$,
$g_{\tilde{C}_i}$ is the genus of the $i$-th component of $\tilde{C}$,
and $C_\sing$ is the set of singular points of $C$.
For example, for the ``conic + tangent,'' formula~(\ref{normalization}) gives
	$$3 = \sum_{i=1}^2 (2-2\cdot 0) + 1 - 2.$$
\end{proof}

\begin{lemma}
\label{smoothpoint}
If $F(x,y,z) \in \F_p[x,y,z]$ is a nonzero homogeneous cubic polynomial
such that $F$ does not factor completely into linear factors over $\Fbar_p$,
then the subscheme $X$ of $\PP^2$ defined by $F=0$
has a smooth $\F_p$-point.
\end{lemma}

\begin{proof}
The polynomial $F$ must be squarefree, since
otherwise $F$ would factor completely.
Hence $X$ is reduced.
If $X$ is a smooth cubic curve,
then it is of genus~1,
and $X(\F_p) \not= \emptyset$ by the Hasse bound.

Otherwise, enumerating possibilities as in Figure~\ref{planecubics}
shows that $X$ is a nodal or cuspidal cubic,
or a union of a line and a conic.
The Galois action on components is trivial,
because when there is more than one,
the components have different degrees.
There is an open subset of $X$
isomorphic over $\F_p$ to $\PP^1$ with at most two geometric points deleted.
But $\#\PP^1(\F_p) \ge 3$, so there remains a smooth $\F_p$-point on $X$.
\end{proof}

\section{The example}
\label{example}

We will carry out the program in Section~\ref{cubicsurface}
with the cubic surface
	$$V: \; 5 x^3 + 9 y^3 + 10 z^3 + 12 w^3 = 0$$
in $\PP^3$.
Cassels and Guy~\cite{casselsguy}
proved that $V$ violates the Hasse principle.
Let $L$ be the line $x+y+z=w=0$.
The intersection $V \cap L$ as a subscheme of $L \isom \PP^1$
with homogeneous coordinates $x,y$ is defined by
	$$5 x^3 + 9 y^3 - 10 (x+y)^3,$$
which has discriminant $242325 = 3^3 \cdot 5^2 \cdot 359 \not= 0$,
so the intersection consists of three distinct geometric points.
This remains true in characteristic $p$,
provided that $p \not\in \{3,5,359\}$.

The projection $V \dashrightarrow \PP^1$ from $L$ is given by
the rational function $u := w/(x+y+z)$ on $V$.
Also, $W$ is the surface in $\PP^3 \times \PP^1$
given by the $((x,y,z,w);(u_0,u_1))$-bihomogeneous equations
\begin{align}
\label{vprime}
	W: \; 5 x^3 + 9 y^3 + 10 z^3 + 12 w^3		&= 0		\\
\nonumber				u_0 w	&= u_1(x+y+z).
\end{align}
The morphism $W \rightarrow \PP^1$
is simply the projection to the second factor,
and the fiber $W_u$ above $u \in \Q = \Aff^1(\Q) \subseteq \PP^1(\Q)$
can also be written as the plane cubic
\begin{equation}
\label{ufiber}
	W_u: \; 5 x^3 + 9 y^3 + 10 z^3 + 12 u^3 (x+y+z)^3 = 0.
\end{equation}

The dehomogenization
	$$h(x,y)=5 x^3 + 9 y^3 + 10 + 12 u^3 (x+y+1)^3,$$
defines an affine open subset in $\Aff^2$ of $W_u$.
Eliminating $x$ and $y$ from the equations
	$$h = \frac{\del h}{\del x} = \frac{\del h}{\del y} = 0$$
shows that this affine variety is singular when $u \in \Qbar$ satisfies
\begin{equation}
\label{singular12}
	2062096 u^{12} + 6065760 u^9 + 4282200 u^6 + 999000 u^3 + 50625 = 0.
\end{equation}
The fiber above $u=0$ is smooth,
so by Lemma~\ref{twelve}, the 12 values of $u$ satisfying~(\ref{singular12})
give the {\em only} points in $\PP^1(\Qbar)$ above which the fiber $W_u$
is singular, and moreover each of these singular fibers is
a nodal cubic.
The polynomial~(\ref{singular12}) is irreducible over $\Q$,
so $W_u$ is smooth for all $u \in \PP^1(\Q)$.

The discriminant of~(\ref{singular12}) is
$2^{146} \cdot 3^{92} \cdot 5^{50} \cdot 359^4$.
Fix a prime $p \not\in \{2,3,5,359\}$,
and a place $\Qbar \dashrightarrow \Fbar_p$.
The 12 singular $u$-values in $\PP^1(\Qbar)$
reduce to 12 distinct singular $u$-values in $\PP^1(\Fbar_p)$
for the family $\Wbar \rightarrow \PP^1$
defined by the two equations~(\ref{vprime}) over $\Fbar_p$.
Moreover, the fiber above $u=0$ is smooth in characteristic $p$.
By Lemma~\ref{twelve},
all the fibers of $\Wbar \rightarrow \PP^1$ in characteristic $p$
are smooth plane cubics or nodal plane cubics.
By Lemma~\ref{smoothpoint} and Hensel's Lemma,
$W_u$ has a $\Q_p$-point for all $u \in \PP^1(\Q_p)$.

\begin{prop}
\label{localpoints}
If $u \in \Q$ satisfies $u \equiv 1 \pmod {p\Z_p}$ for $p \in \{2,3,5\}$
and $u \in \Z_{359}$,
then the fiber $W_u$ has a $\Q_p$-point for all completions $\Q_p$,
$p \le \infty$.
\end{prop}

\begin{proof}
Existence of real points is automatic,
since $W_u$ is a plane curve of odd degree.
Existence of $\Q_p$-points for $p \not\in \{2,3,5,359\}$
was proved just above the statement of Proposition~\ref{localpoints}.

Consider $p=359$.
A Gr\"obner basis calculation shows that
there do not exist
$a_1$, $a_2$, $b_1$, $b_2$, $c_1$, $c_2$, $\overline{u} \in \Fbar_{359}$
such that
\begin{equation}
\label{ubarfiber}
	5 x^3 + 9 y^3 + 10 z^3 + 12 \overline{u}^3 (x+y+z)^3
\end{equation}
and
	$$(5x + a_1 y + a_2 z)(x + b_1 y + b_2 z)(x + c_1 y + c_2 z)$$
are identical.
Hence Lemma~\ref{smoothpoint} applies to show that
for any $\overline{u} \in \F_{359}$,
the plane cubic
defined by~(\ref{ubarfiber}) over $\F_{359}$ has a smooth $\F_{359}$-point,
and Hensel's Lemma implies that $W_u$ has a $\Q_{359}$-point
at least when $u \in \Z_{359}$.

When $u \equiv 1 \pmod{5\Z_5}$,
the curve reduced modulo~5,
	$$\Wbar_{\overline{u}}: \; 4 y^3 + 2(x+y+z)^3 = 0,$$
consists of three lines through $P:=(1:0:-1) \in \PP^2(\F_5)$,
so it does not satisfy the conditions of Lemma~\ref{smoothpoint},
but one of the lines, namely $y=-2(x+y+z)$, is defined over $\F_5$,
and every $\F_5$-point on this line except $P$
is smooth on $\Wbar_{\overline{u}}$.
Hence $W_u$ has a $\Q_5$-point.

The same argument shows that $W_u$ has a $\Q_2$-point
whenever $u \equiv 1 \pmod{2\Z_2}$,
since the curve reduced modulo~2 is $x^3 + y^3 = 0$,
which contains $x+y=0$.

Finally, when $u \equiv 1 \pmod{3\Z_3}$,
the point $(1:2:1)$ satisfies the equation~(\ref{ufiber})
modulo~$3^2$, and Hensel's Lemma gives a point $(x_0:2:1) \in W_u(\Q_3)$
with $x_0 \equiv 1 \pmod{3\Z_3}$.
This completes the proof.
\end{proof}

We now seek a non-constant rational function $\PP^1 \rightarrow \PP^1$
that maps $\PP^1(\Q_p)$ into $1+p\Z_p$ for $p \in \{2,3,5\}$
and into $\Z_{359}$ for $p=359$.
For $p \in \{2,3,5\}$,
the rational function $v=t^4$ maps $t \in \PP^1(\F_p)$ into $\{0,1,\infty\}$,
and
	$$u=\frac{v^3-v-1}{v^3-v^2-1}$$
maps any $v \in \{0,1,\infty\}$ to $1$.
Moreover $u \not= \infty$ for $v \in \PP^1(\F_{359})$.
Therefore, for $t \in \PP^1(\Q)$,
the value of
	$$u = \frac{t^{12}-t^4-1}{t^{12}-t^8-1}$$
satisfies the local conditions in Proposition~\ref{localpoints}.

Substituting into~(\ref{ufiber}),
we see that
\begin{equation}
\label{tfiber}
	X_t: \; 5 x^3 + 9 y^3 + 10 z^3
	+ 12 \left( \frac{t^{12}-t^4-1}{t^{12}-t^8-1} \right)^3 (x+y+z)^3
	= 0
\end{equation}
has $\Q_p$-points for all $p \le \infty$.
On the other hand, $X_t(\Q) = \emptyset$, because $V(\Q) = \emptyset$.
Finally, the existence of nodal fibers
in the family implies as in~\cite{colliotpoonen}
that the $j$-invariant of the family has poles, and hence is non-constant.

\section{The Jacobians}
\label{jacobians}

For $t \in \Q$,
let $E_t$ denote the Jacobian of $X_t$.
Salmon in his work on invariants of a plane cubic
developed formulas which were later shown in~\cite{genus1}
to be coefficients of a Weierstrass model of the Jacobian.
We used a GP-PARI implementation of these by Fernando Rodriguez-Villegas,
available electronically at
ftp://www.ma.utexas.edu/pub/villegas/gp/inv-cubic.gp
to show that our $E_t$ has a Weierstrass model $y^2 = x^3 + A x + B$
where
\begin{align*}
	A &= 145800 (t^{12}-t^4-1)^3 (t^{12}-t^8-1) \\
\intertext{and}
	B &=  -6129675 t^{72} + 48660750 t^{68} - 72845325 t^{64} - 43776450 t^{60} - 52032375 t^{56} \\
	&\qquad + 392384250 t^{52} - 12636000 t^{48} - 198105750 t^{44} - 604705500 t^{40} + 229027500 t^{36} \\
	&\qquad + 387518175 t^{32} + 384183000 t^{28} - 242872425 t^{24} - 227776050 t^{20} - 110899125 t^{16} \\
	&\qquad + 76387050 t^{12} + 41735250 t^8 + 11882700 t^4 - 6129675.
\end{align*}
Because the non-existence of rational points on $V$
is explained by a Brauer-Manin obstruction,
Section~3.5 and in particular Proposition~3.5 of~\cite{colliotpoonen}
show that there exists a {\em second family}
of genus~1 curves $Y_t$ with the same Jacobians
such that the Cassels-Tate pairing satisfies
$\langle X_t, Y_t \rangle = 1/3$ for all $t \in \Q$.
In particular, for all $t \in \Q$,
the Shafarevich-Tate group $\Sha(E_t)$
contains a subgroup isomorphic to $\Z/3 \times \Z/3$.

\section*{Acknowledgements}

I thank Ahmed Abbes for explaining the formula~(\ref{ecformula}) to me.
I thank also Bill McCallum and Fernando Rodriguez-Villegas,
for providing Salmon's formulas for invariants of plane cubics
in electronic form.
The calculations for this paper were mostly done using Mathematica
and GP-PARI on a Sun Ultra~2.
The values of $A$ and $B$ in Section~\ref{jacobians}
and their transcription into LaTeX
were checked by pasting the LaTeX formulas into Mathematica,
plugging them into the formulas for the $j$-invariant from GP-PARI,
and comparing the result
against the $j$-invariant of $X_t$
as computed directly by Mark van Hoeij's Maple package ``IntBasis'' at 
http://www.math.fsu.edu/~hoeij/compalg/IntBasis/index.html
for a few values of $t \in \Q$.

\end{document}